\documentstyle[12pt]{article}
\newcommand{\parf}{\subsection}
\newcommand{\punk}{\subsubsection}
\font\bfit=cmbxti10 scaled \magstep1
\expandafter\chardef\csname pre amsym.def 
at\endcsname=\the\catcode`\@
\catcode`\@=11
\def\newsymbol#1#2#3#4#5{\let\next@\relax
\ifnum#2=\@ne\let\next@\msafam@\else
 \ifnum#2=\tw@\let\next@\msbfam@\fi\fi
 \mathchardef#1="#3\next@#4#5}
\def\hexnumber@#1{\ifcase#1 0\or 1\or 2\or 3\or 4\or 5\or 6\or 
7\or 8\or 9\or A\or B\or C\or D\or E\or F\fi}
\newfam\msafam
\newfam\msbfam
\edef\msafam@{\hexnumber@\msafam}
\edef\msbfam@{\hexnumber@\msbfam}
\def\Bbb#1{\fam\msbfam\relax#1}
\catcode`\@=\csname pre amsym.def at\endcsname
\newsymbol\subsetneqq  2324
\newsymbol\restriction 1316
\newsymbol\complement  107B
\newsymbol\preceQ      1334
\newsymbol\succeQ      133C
\newsymbol\squarE      1003
\font\twlmsa=msam10 scaled \magstep1
\font\tenmsa=msam10 
\font\egtmsa=msam8
\def\xiimsa{\textfont\msafam=\twlmsa\scriptfont\msafam=\egtmsa}
\def\xmsa  {\textfont\msafam=\tenmsa\scriptfont\msafam=\egtmsa}
\font\twlmsb=msbm10 scaled \magstep1
\font\tenmsb=msbm10
\font\egtmsb=msbm8
\def\xiimsb{\textfont\msbfam=\twlmsb\scriptfont\msbfam=\egtmsb}
\def\xmsb  {\textfont\msbfam=\tenmsb\scriptfont\msbfam=\egtmsb}

\font\twleufm=eufm10 scaled \magstep1

\font\egteufm=eufm8
\newfam\eufmfam
\newcommand{\xiieufm}{\textfont\eufmfam=\twleufm
    \scriptfont\eufmfam=\egteufm\def\fraK{\fam\eufmfam\twleufm}}

\font\twleusm=eusm10 scaled \magstep1
\font\teneusm=eusm10
\font\egteusm=eusm8
\newfam\eusmfam
\newcommand{\xiieusm}{\textfont\eusmfam=\twleusm
    \scriptfont\eusmfam=\egteusm\def\skr{\fam\eusmfam\twleusm}}
\newcommand{\xeusm}{\textfont\eusmfam=\teneusm
\scriptfont\eusmfam=\egteusm\def\skr{\fam\eusmfam\teneusm}}

\font\ninss=cmss9
\font\egtss=cmss8

\newcommand{\xiiss} {\scriptfont\sffam=\ninss}
\newcommand{\xss}   {\scriptfont\sffam=\egtss}

\def\addto#1#2{
\ifx\zzone\undefined\let\zzone=#1\def#1{\zzone#2}\else
\ifx\zztwo\undefined\let\zztwo=#1\def#1{\zztwo#2}\else
\fi\fi
}
\addto\normalsize   {\xiimsb\xiimsa\xiieufm\xiieusm\xiiss} 
\addto\footnotesize {\xmsb\xmsa\xeusm\xss} 
\newtheorem{theorem}             {Theorem}
\newtheorem{corollary}  [theorem]{Corollary}
\newtheorem{definition} [theorem]{Definition}
\newtheorem{lemma}      [theorem]{Lemma}
\newtheorem{proposition}[theorem]{Proposition}
\newtheorem{aq}    {{\bfit Acknowledgements}}  
\newtheorem{prF}   {Proof}                     
\newcommand{\thsp}{\hspace{0.1ex}}
\newcommand{\bte} {\begin{theorem}}
\newcommand{\ete} {\end{theorem}}
\newcommand{\bcor}{\begin{corollary}\thsp}
\newcommand{\ecor}{\end{corollary}}
\newcommand{\bdf} {\begin{definition}\rm\thsp}
\newcommand{\edf} {\end{definition}}
\newcommand{\ble} {\begin{lemma}\thsp}
\newcommand{\ele} {\end{lemma}}
\newcommand{\baq} {\begin{aq}\rm\thsp}
\newcommand{\eaq} {\end{aq}}
\newcommand{\bpf} {\begin{prF}\rm} 
\newcommand{\epf} {\qeD\end{prF}} 
\newcommand{\bpro}{\begin{proposition}\thsp} 
\newcommand{\epro}{\end{proposition}} 
\newcommand{\qeD}{\hspace{0mm}\hfill$\mtho\squarE$}
\newcommand{\ben}{\begin{enumerate}}
\newcommand{\een}{\end{enumerate}}
\newcommand{\bay}{\begin{array}}
\newcommand{\eay}{\end{array}}
\newcommand{\dosp}{\hspace{0.5ex}}
\newcommand{\disp}{\hspace{0.35ex}}
\newcommand{\ea}  {\hbox{\it e.\disp a.\/}}
\newcommand{\ie}  {\hbox{\it i.\disp e.\/}}
\newcommand{\etc} {{\it etc.\/}}
\newcommand{\eg}  {\hbox{\it e.\disp g.\/}}
\newcommand{\po}  {\hbox{p.\dosp o.}}
\newcommand{\pqo} {\hbox{p.\dosp q.-o.}}
\newcommand{\lo}  {\hbox{l.\dosp o.}}

\newcommand{\hop} {\hbox{h.\dosp o.\dosp p.}}

\newcommand{\mek}  {\preceQ}
   
\newcommand{\mso} {\mathbin{<_0}} 
\newcommand{\meo} {\mathbin{\leq_0}}
\newcommand{\msl} {\mathbin{<_{\rbox{lex}}}} 
\newcommand{\mel} {\mathbin{\leq_{\rbox{lex}}}}
\newcommand{\eee} {\approx}
\newcommand{\eqf} {\equiv} 
\newcommand{\gc} {{\fraK c}}
\newcommand{\gP} {{\fraK P}}
\newcommand{\gp} {{\fraK p}}

\newcommand{\pri}{{\rbox{\tt pr}}_1\hspace{1pt}}
\newcommand{\prt}{{\rbox{\tt pr}}_2\hspace{1pt}}
\newcommand{\dvo}[1]{{\Bbb#1}}
\newcommand{\dP}    {\dvo P}

\newcommand{\dZ}    {\dvo Z}
\newcommand{\dpp}{\dP^+_2}
\newcommand{\dpm}{\dP^-_2}
\newcommand{\dpt}{\dP^2_{\eqf}}
\newcommand{\skrsp}{\hspace{0ex}}
\newcommand{\skri}[1]{{\skrsp \skr{#1}\skrsp}}
\newcommand{\cF} {{\skri F}}
\newcommand{\cN} {{\skri N}}
\newcommand{\cO} {{\skri O}}
\newcommand{\tfu}{\cF}
\newcommand{\itla}{\item\label}
\newcommand{\rbox} [1]{{\rm{#1}}}        
\newcommand{\bbox} [1]{{\bf{#1}}}

\newcommand{\kV} {\bbox{V}}
\newcommand{\kvp}{\kV^+}
\newcommand{\al}  {\alpha} 
\newcommand{\ba}  {\beta} 

\newcommand{\Ups} {\Upsilon}
\newcommand{\kpa} {\kappa}
\newcommand{\om}  {\omega} 
\newcommand{\omi} {\om_1} 
\newcommand{\lom} {^{<\om}} 
 
\newcommand{\omck}{\om_1^{\rbox{CK}}}
\newcommand{\iSg}{{\mathchar"7106}}
\newcommand{\iPi}{{\mathchar"7105}}
\newcommand{\iDa}{{\mathchar"7101}}
\newcommand{\is}[2]{\iSg^{#1}_{#2}}
\newcommand{\ip}[2]{\iPi^{#1}_{#2}}
\newcommand{\id}[2]{\iDa^{#1}_{#2}}
\newcommand{\ti}  {\times}
\newcommand{\res} {{\hspace{0.2ex}\restriction\hspace{0.2ex}}}
\newcommand{\pone}{\hspace{-0.4ex}+\hspace{-0.4ex}1}
\newcommand{\sq}  {\subseteq}
\newcommand{\cj}  {\mathbin{\hspace{0.2ex}\wedge\hspace{0.2ex}}}
\newcommand{\orr} {\mathbin{\textstyle\bigvee}}
\newcommand{\imp} {\mathbin{\,\Longrightarrow\,}}
\newcommand{\eqv} {\mathbin{\,\Longleftrightarrow\,}}
\newcommand{\lra} {\longrightarrow} 
\newcommand{\ans} [1] {\{\hspace{0.1ex}#1\hspace{0.1ex}\}}
\newcommand{\strk}[2] {{\ang{#1\hspace{1pt};\hspace{1pt}#2}}}

\newcommand{\ang} [1] {\langle #1\rangle}
\newcommand{\ima}{\mathbin{\hbox{\rm''}}}
\newcommand{\we}  {{\mathbin{\hspace*{0.2ex}^\wedge}}}
\newcommand{\sus} {{\exists\,}}
\newcommand{\kaz} {{\forall\,}}
\newcommand{\dm}  {$$}
\newcommand{\sneq}{\subsetneqq}
\newcommand{\emps}{\emptyset}
\newcommand{\relf}{\sf}
\newcommand{\qE}  {\mathbin{\relf E}}
\newcommand{\Eo}  {\mathbin{{\qE}_0}}
\newcommand{\nEo} {\mathbin{\not{{\hspace{-0.4ex}\qE}}_0}}
\newcommand{\Xin}  {X_\infty}
\newcommand{\mtho}{\mathsurround=0mm}
\newcommand{\msur}{\hspace*{-1\mathsurround}}

\newcommand{\noi}{\noindent}
\newcommand{\vom}{\vspace{1mm}}
\newcommand{\vtm}{\vspace{2mm}}
\newcommand{\dd}[2]{\hbox{$\mtho\hspace{0.2ex}{#1}$-#2}}

\begin{document}

\normalsize

\title{When a partial Borel order is linearizable}

\author{Vladimir Kanovei
\thanks{\ Moscow Transport Engineering Institute}
\thanks{\ {\tt kanovei@mech.math.msu.su} \ and \ 
{\tt kanovei@math.uni-wuppertal.de}
}
\thanks{\ This paper was accomplished during my visit to Caltech 
in April 1997. I thank Caltech for the support and A. S. Kechris 
and J. Zapletal for useful information and interesting discussions 
relevant to the topic of this paper during the visit.}
}
\date{April 1997} 
\maketitle
\normalsize

\begin{abstract}
We prove the following classification theorem of the 
``Glimm -- Effros'' type for Borel order relations: a Borel 
partial order on the reals either is Borel linearizable or 
includes a copy of a certain Borel partial order $\meo$ which 
is not Borel linearizable. 
\end{abstract}

\newpage

\noi
\parf*{Notation}
\label{not}

A binary relation $\mek$ on a set $X$ is a 
{\it partial quasi-order\/}, or {\it\pqo\/} in brief, on $X,$ iff 
$x\mek y\cj y\mek z\imp x\mek z,$ and $x\mek x$ for any $x\in X.$ 
In this case, $\eee$ is the associated equivalence relation, 
\ie\ $x\eee y$ iff $x\mek y\cj y\mek x$. 

If in addition $x\eee x\imp x=x$ for any $x$ then $\mek$ is 
a {\it partial order\/}, or {\it\po\/}, so that, say, 
forcing relations are \pqo's, but, generally speaking, not \po's 
in this notational system. 

A 
\po is {\it linear\/} 
(\lo) iff 
we have $x\mek y\orr y\mek x$ for all $x,\,y\in X$.

Let $\mek$ and $\mek'$ be \pqo's on resp.\ $X$ and $X'.$ A map 
$h:X \lra X'$ will be called {\it half order preserving\/}, 
or {\it\hop\/}, iff 
$
x\mek y\imp h(x)\mek' h(y),
$.

\bdf
A Borel \pqo\ $\strk X\mek$ is {\it Borel linearizable\/} iff 
there is a Borel \lo\ $\strk{X'}{\mek'}$ and a Borel \hop\ map 
$h:X\lra X'$ (called: {\it linearization map\/}) satisfying 
$x\eee y\eqv h(x)=h(y)$.~\footnote
{\rm\ The equivalence cannot be dropped as otherwise a 
one-element set $X'$ works in any case.}
\edf

\parf*{Introduction}
\label{intr}

Harrington, Marker, and Shelah \cite{hms} proved several theorems 
on Borel partial orders, mainly concerning {\it thin\/} \pqo's, 
\ie\ those which do not admit uncountable pairwise incomparable 
subsets. In particular they demonstrated that any such a Borel 
\pqo\ is Borel linearizable, and moreover the corresponding \lo\ 
$\strk{X'}{\mek'}$ can be chosen so that $X'\sq 2^\al$ for some 
$\al<\omi$ while ${\mek'}={{\mel}\res X},$ where $\mel$ is the 
lexicographical order. 

As elementary examples show that the thinness is not a necessary 
condition for the Borel linearization, this result leaves open 
the problem of linearization of non-thin Borel \po's. Harrington 
\ea\ write in \cite{hms} that ``there is little to say about 
nonthin orderings'', although there are many interesting among 
them like the {\it dominance\/} order on $\om^\om.$ 

Our main result will say that not all Borel \pqo's are Borel 
linearizable, and there exists a {\it minimal\/} in certain 
sense among them. 

\bdf
Let $a,\,b\in 2^\om.$ We define $a\meo b$ iff either $a=b$ or 
$a\Eo b$~\footnote
{\rm\ That is $a(k)=b(k)$ for all but finite $k$ --- the 
{\it Vitali\/} equivalence relation on $2^\om.$} 
and $a(k_0)<b(k_0)$ where $k_0$ is the largest $k$ such that 
$a(k)\not=b(k)$.~\footnote
{\rm\ If one enlarges $\mso$ so that, in addition, 
$a\mso b$ whenever $a,\,b\in 2^\om$ are such that $a(k)=1$ and 
$b(k)=0$ for all but finite $k$ then the enlarged relation can 
be induced by a Borel action of $\dZ$ on $2^\om,$ such that 
$a\mso b$ iff $a=zb$ for some $z\in \dZ,\msur$ $z>0$.}

\edf
The relation $\meo$ is a Borel \pqo\ on $2^\om$ which orders 
every \dd\Eo class similarly to the integers $\dZ$ (except for 
the class $[\om\ti\ans{0}]_{\Eo}$ ordered as 
$\om$ and the 
\pagebreak[1] 
class $[\om\ti\ans{1}]_{\Eo}$ ordered as $\om^\ast$) 
but leaves any two \dd\Eo inequivalent reals incomparable.

\bte
\label{hb}
{\rm(The main result.)} \ 
Suppose that\/ ${\mek}$ is a Borel \pqo\ on\/ $\cN=\om^\om.$ 
Then exactly one of the following two conditions is satisfied$:$
\ben
\itemsep=1mm
\def\theenumi{(\Roman{enumi})}
\def\labelenumi{{\rm\theenumi}}
\itla{1B}\msur
$\mek$ is Borel linearizable, moreover there exist an ordinal\/ 
$\al<\omi$ and a Borel linearization map\/ 
$h:\strk\cN\mek\lra \strk{2^\al}\mel$.

\itla{2B} 
there exists a continuous\/ $1-1$ map\/ $F:2^\om\lra\cN$ 
such that we have\/ $a\meo b\imp F(a)\mek F(b)$ while\/ 
${a\nEo b}$ implies that $F(a)$ and\/ $F(b)$ are\/ \dd\mek 
incomparable.~\footnote
{\rm\ Then $F$ associates a chain 
$\ans{F(b):b\Eo a}$ in $\strk\cN\mek$ to each \dd\Eo class 
$[a]_{\Eo}$ so that any two different chains do not contain 
\dd\mek comparable elements: let us call them 
{\it fully incomparable\/} chains. Thus \ref{2B} essentially says 
that $\mek$ admits an effectively ``big'' Borel family of fully 
incomparable chains, which is therefore necessary and 
sufficient for $\mek$ to be {\it not\/} Borel 
linearizable.}
\een
\ete
The theorem resembles the case of Borel equivalence relations 
where a necessary and sufficient condition for a Borel 
equivalence relation $\qE$ to be {\it smooth\/} is that $\Eo$ 
(which is not smooth) does not continuously embed in $\qE$ 
(Harrington, Kechris, Louveau \cite{hkl}). 
($\meo$ itself is {\it not\/} Borel linearizable.) 

The proof is essentially a combination of ideas 
and technique in \cite{hkl,hms}.

\parf{Incompatibility}
\label{inc}

Let us first prove that \ref{1B} and \ref{2B} are incompatible. 

{\it Suppose otherwise\/}. 
The superposition of the maps $F$ and $h$ is then a Borel 
\hop\ map $\phi:\strk{2^\om}{\meo}\lra \strk{2^\al}{\mel}$ 
satisfying the following: $\phi(a)=\phi(b)$ implies that 
$a\Eo b,$ \ie\ $a$ and $b$ are \dd\meo comparable. 

Therefore, as any \dd\Eo class is \dd\meo ordered similarly 
to $\dZ,\msur$ $\om,$ or $\om^\ast,$ the \dd\phi image 
$X_a=\phi\ima{[a]_{\Eo}}$ of the \dd{\Eo}class of any 
$a\in 2^\om$ is \dd\mel ordered similarly to a subset of $\dZ.$ 
If $X_a=\ans{x_a}$ is a singleton then put $\psi(a)=x_a$.  

Assume now that $X_a$ contains at least two points. In this case 
we can effectively pick an element in $X_a\,!$ Indeed there is a 
maximal sequence $u\in 2^{<\al}$ such that $u\subset x$ for each 
$x\in X_a.$ Then the set 
$X_a^{\rm left}=\ans{x\in X:u\we 0\subset x}$ 
contains a \dd\mel largest element, which we denote by $\psi(a)$. 

To conclude $\psi$ is a Borel reduction of $\Eo$ to the equality 
on $2^\al,$ \ie\ $a\Eo b$ iff $\psi(a)=\psi(b),$ which is 
impossible.

\parf{The dichotomy}
\label d

As usual it will be assumed that the \pqo\ $\mek$ of 
Theorem~\ref{hb} is a $\id11$ 
relation. Let $\eee$ denote the associated equivalence.

Following \cite{hms} let, for $\al<\omck,$ $\tfu_\al$ to 
be the family of all \hop\ $\id11$ functions 
$f:\strk\cN\mek\lra \strk{2^\al}\mel.$ 
Then $\tfu=\bigcup_{\al<\omck}\tfu_\al$ is a (countable) $\ip11$ 
set, in a suitable coding system for functions of this type.
(See \cite{hms} on details.)

Define, for $x,\,y\in\cN,$ $x\eqf y$ iff $f(x)=f(y)$ for any 
$f\in \tfu.$ 

\ble
\label{eqfs}
{\rm(See \cite{hms}.)} \ 
$\eqf$ is a\/ $\is11$ equivalence relation including\/ $\eee$.
\ele
\bpf
As $\mek$ is $\id11,$ one gets by a rather standard argument a set 
$N\sq\om$ and a function $f_n\in\tfu$ for any $n\in N$ so that 
$\tfu=\ans{f_n:n\in N}$ and the relations  
$n\in N\cj f_n(x)\mel f_n(y)$ and $n\in N\cj f_n(x)\msl f_n(y)$
are presentable in the form 
$n\in N\cj\cO(x,y)$ and $n\in N\cj\cO'(x,y)$ 
where $\cO,\,\cO'$ are $\is11$ relations. Now 
$x\eqf y$ iff $\kaz n\:(n\in N\imp f_n(x)=f_n(y)),$ as required. 
\epf

{\bfit Case 1\/\bf:\/} $\eqf$ coincides with $\eee$.\vtm

Let us show how this implies \ref{1B} of Theorem~\ref{hb}. The 
set
\dm
P=\ans{\ang{x,y,n}:x\not\eee y\cj f_n(x)\not=f_n(y)}.
\dm
is $\ip11$ and, by the assumption of Case 1, its projection on 
$x,\,y$ coincides with the complement of $\eee.$ Let $Q\sq P$ be 
a $\ip11$ set uniformizing $P$ in the sense $\cN^2\ti \om.$ Then 
$Q$ is $\id11$ because  
$
Q(x,y,n)\eqv x\not\eee y\cj\kaz n'\not=n\:\neg\:Q(x,y,n').
$
It follows that $N'=\ans{n:\sus x,y\:Q(x,y,n)}\sq N$ is $\is11.$ 
Therefore there is a $\id11$ set $M$ such that $N'\sq M\sq N.$ 
\footnote
{\ Harrington \ea~\cite{hms} use a general reflection theorem to 
get such a set, but a more elementary reasoning sometimes has 
advantage.} 

Consider a $\id11$ enumeration $M=\ans{n_l:l\in\om}.$ For any $l,$ 
$f_{n_l}\in\tfu_\al$ for some ordinal $\al=\al_l<\omck.$ Another 
standard argument (see \cite{hms}) shows that in this case (\eg\ 
when $M\sq N$ is $\id11$) the ordinals $\al_l$ are bounded by 
some $\al<\omck.$ It follows that the function
$
h(x)=f_{n_0}(x)\we f_{n_1}(x)\we f_{n_2}(x)\we ... \we 
f_{n_l}(x)\we ...
$
belongs to some $\tfu_\ba,\msur$ $\ba\leq\al\cdot \om.$ On the 
other hand, by the construction we have $x\eee y\eqv h(x)=h(y)$ 
so $h$ satisfies \ref{1B} of Theorem~\ref{hb}.\vtm

{\bfit Case 2\/\bf:\/} ${\eee}\sneq{\eqf}$. 
Assuming this we work towards \ref{2B} of Theorem~\ref{hb}. 

\parf{The domain of singularity}
\label{sing}

By the assumption the $\is11$ set 
$A=\ans{x:\sus y\:(x\eee y\cj x\not\eqf y)}$ is non-empty. 

Define $X\eqf Y$ iff $\kaz x\in X\:\sus y\in Y\:x\eqf y$ and 
{\it vice versa\/}. 

\bpro
\label p
Let\/ $X,\,Y\sq A$ be non-empty\/ $\is11$ sets satisfying\/ 
$X\eqf Y.$ Then 
\dm
P_+=\ans{\ang{x,y}\in X\ti Y:x\eqf y\cj x\mek y}
\hspace{1pt},\hspace{3.5pt}
P_-=\ans{\ang{x,y}\in X\ti Y:x\eqf y\cj x\not\mek y}
\dm
are non-empty\/ $\is11$ sets, their projections\/ 
$\pri P^+$ and\/ $\pri P^-$ 
are \dd{\is11}dense in\/ $X$~\footnote
{\rm\ That is intersect any non-empty $\is11$ set $X'\sq X$.}
while the projections\/ $\prt P^+$ and\/ $\prt P^-$ are\/ 
\dd{\is11}dense in\/ $Y$.~\footnote
{\ For a set $P\sq\cN^2,$ $\pri P$ and $\prt P$ have the 
obvious meaning of the projections on the resp.\ 1st and 2nd copy 
of $\cN$.}
\epro
\bpf
The density easily follows from the non-emptiness, so let us 
concentrate on the latter. {\it Prove that\/ $P_+\not=\emps$.}

Suppose on the contrary that $P_+=\emps.$ Then there is a single 
function $f\in\tfu$ such that the set 
$
\ans{\ang{x,y}\in X\ti Y:f(x)=f(y)\cj x\mek y}
$
is empty. (See the reasoning in Case 1 of Section~\ref{d}.) Define 
\pagebreak[0]
\dm
\Xin=\ans{x:\kaz y\in Y\:(f(x)=f(y)\imp x\not\mek y)},
\dm
so that $\Xin$ is a $\ip11$ set and $X\sq\Xin$ but 
$Y\cap\Xin=\emps.$ Using Separation we can easily define a 
sequence of sets
\dm
X=X_0\sq U_0\sq X_1\sq U_1\sq\dots\sq X_n\sq U_n\sq\dots\sq\Xin
\dm
so that $U_n=\ans{x':\sus x\in X_n\:(f(x)=f(x')\cj x\mek x')}$ 
while $X_{n+1}\in\id11$ for all $n.$ (Note that if $X_n\sq\Xin$ 
and $U_n$ is defined as indicated then $U_n\sq\Xin$ too.) Moreover 
a proper execution of the construction~\footnote
{\ We refer to the proof of an ``invariant'' effective Separation 
theorem in Harrington, Kechris, Louveau~\cite{hkl}, which includes 
a similar construction.}
allows to get the final set $U=\bigcup_n U_n=\bigcup_n X_n$ in 
$\id11.$ Note that $X\sq U$ but $Y\cap U=\emps$ as $U\sq\Xin$.

Now put $f'(x)=f(x)\we 1$ whenever $x\in U$ and $f'(x)=f(x)\we 0$ 
otherwise. {\it We assert that\/ $f'\in\tfu.$\/} Indeed suppose 
that $x'\mek y'$ and prove $f'(x')\mel f'(y').$ It can be 
assumed that $f(x')=f(y').$ It remains to check 
$x'\in U\imp y'\in U,$ which easily follows from the definition 
of sets $U_n.$ Thus $f'\in\tfu$.

However clearly $f'(x)\not=f'(y),$ hence $x\not\eqf y,$ whenever 
$x\in X$ and $y\in Y$ 
which is a contradiction with the assumption that $X\eqf Y$. 

Now {\it prove that\/ $P_-\not=\emps.$\/} Consider first the 
case $X=Y.$ Suppose on the contrary that $P_-=\emps.$ Then, as 
above, there is a single function $f\in\tfu$ such that the set 
$
\ans{\ang{x,y}\in X^2:f(x)=f(y)\cj x\not\mek y}
$
is empty, so that $\eqf$ and $\eee$ coincide on $X.$ Our plan is 
to find functions $f',\,f''\in\tfu$ such that the sets
\dm
\bay{rcl}
Q'&=&\ans{\ang{x,y}\in X\ti\cN:f'(x)=f'(y)\cj y\not\mek x}\\[0.5em]
Q''&=&\ans{\ang{x,y}\in X\ti\cN:f''(x)=f''(y)\cj x\not\mek y}
\eay
\dm
are empty; then 
$Q=\ans{\ang{x,y}\in X\ti\cN:x\eqf y\cj y\not\eee x}=\emps,$ 
which a contradiction with $\emps\not=X\sq A$. 

Let us find $f';$ the other case is similar. Define 
\dm
\Xin=\ans{x:\kaz x'\in X\:(f(x)=f(x')\imp x\mek x')},
\dm
so that $\Xin$ is $\ip11$ and $X\sq\Xin.$ As above there is a 
sequence of sets
\dm
X=X_0\sq U_0\sq X_1\sq U_1\sq\dots\sq X_n\sq U_n\sq\dots\sq\Xin
\dm
such that $U_n=\ans{u:\sus x\in X_n\:(f(x)=f(u)\cj u\mek x)}$ 
while $X_{n+1}\in\id11$ for all $n$ and the final set 
$U=\bigcup_n U_n=\bigcup_n X_n$ belongs to $\id11$.

Set $f'(x)=f(x)\we 0$ whenever $x\in U$ and $f'(x)=f(x)\we 1$ 
otherwise. Then $f'\in\tfu.$ Prove that $f'$ witnesses that 
$Q'=\emps.$ Consider any $x\in X$ and $y\in\cN$ such that 
$f'(x)=f'(y).$ Then in particular $f(x)=f(y)$ and 
$x\in U\eqv y\in U,$ so that $y\in U$ because we know that 
$x\in X\sq U.$ Thus $U\in\Xin,$ so by definition $y\mek x$ 
as required.

Finally prove $P_-\not=\emps$ in the general case. By the 
result in the particular case, the $\is11$ set
$
P'=\ans{\ang{x,x'}\in X^2:x\eqf x'\cj x\not\mek x'}
$
is non-empty. Let $X'=\ans{x'\in X:\sus x\:P'(x,x')}$ and 
$Y'=\ans{y\in Y:\sus x'\in X'\:(x'\eqf y)},$ so that $X'$ and 
$Y'$ are $\is11$ sets satisfying $X'\eqf Y'.$ By the result for 
$P_+$ there exist $x'\in X'$ and $y\in Y'$ satisfying $x'\eqf y$ 
and $y\mek x'.$ Now there exists $x\in X$ such that $x\eqf x'$ 
and $x\not\mek x'.$ Then $x\eqf y$ and $x\not\mek y,$ as required.
\epf

\parf{The forcing notions involved}
\label f

Our further strategy will be the following. We shall 
define a generic extension of the universe $\kV$ (where 
Theorem~\ref{hb} is being proved) in which there exists a 
function $F$ which witnesses \ref{2B} of Theorem~\ref{hb}. 
However as the existence of such a function is a $\is12$ 
statement, we obtain the result for $\kV$ by Shoenfield.
\footnote
{\ In fact the proof can be conducted without any use of 
metamathematics, as in Harrington \ea\ \cite{hkl}, but at the 
cost of longer reasoning.}

\bdf 
$\dP$ is the collection of all non-empty $\is11$ sets $X\sq A.$ 
\edf
It is a standard fact that $\dP$ (the {\it Gandy forcing\/}) forces 
a real which is the only real which belongs to every set in the 
generic set $G\sq \dP.$ (We identify $\is11$ sets in the universe 
$\kV$ with their copies in the extension.)

\bdf
$\dpp$ is the collection of all non-empty $\is11$ sets 
$P\sq A^2$ such that $P(x,y)\imp x\eqf y\cj x\mek y.$  
$\dpm$ is defined similarly but with $x\not\mek y$ instead. 
\edf
Both of them 
force a pair of reals $\ang{x,y}\in A^2$ satisfying resp.\ 
$x\mek y$ and $x\not\mek y.$ ($\dpp$ and $\dpm$ are non-empty 
forcing notions by Proposition~\ref p.) 

\bdf
$\dpt$ is the collection of all sets of the form 
$\Ups=X\ti Y$ where $X,\,Y$ are sets in $\dP$ satisfying 
$X\eqf Y$.
\edf

\ble
\label{pp}
$\dpt$ forces a pair of reals\/ $\ang{x,y}$ such that\/ 
$x\not\mek y$.
\ele
\bpf
Suppose that on the contrary a condition $\Ups_0=X_0\ti Y_0$ 
in $\dpt$ forces $x\mek y.$ Consider a more complicated forcing 
\pagebreak[1]
$\gP$ which consists of forcing conditions of the form 
$\gp=\ang{\Ups,P,\Ups',Q},$ where $\Ups=X\ti Y$ and 
$\Ups'=X'\ti Y'$ belong to $\dpt,\msur$ $P\in\dpp,\msur$ 
$P\sq Y\ti X',\msur$ $Q\in\dpm,\msur$ $Q\sq X\ti Y',$ 
and the sets $\pri P\sq Y,\msur$ $\prt P\sq X',\msur$  
$\pri Q\sq X$ and $\prt Q\sq Y'$ are \dd{\is11}dense in resp.\ 
$Y,\,X',\,X,\,Y'.$ 

For instance setting 
$P_0=\ans{\ang{y,x'}\in Y_0\ti X_0:y\eqf x'\cj y\mek x'}$ and 
$Q_0=\ans{\ang{x,y'}\in X_0\ti Y_0:x\eqf y'\cj x\not\mek y'}$ 
we get a condition $\gp_0=\ang{\Ups_0,P_0,\Ups_0,Q_0}\in\gP$ by 
Proposition~\ref p.

It is the principal fact that if $\gp=\ang{\Ups,P,\Ups',Q}\in\gP$ 
and we strengthen one of the components within the corresponding 
forcing notion then this can be appropriately reflected in the 
other components. To be concrete assume that, for instance, 
$P^\ast\in\dpp,\msur$ $P^\ast\sq P,$ and find a condition 
$\gp_1=\ang{\Ups_1,P_1,\Ups'_1,Q_1}\in\gP$ satisfying 
$\Ups_1\sq\Ups,\msur$ $\Ups'_1\sq\Ups',\msur$ $P_1\sq P^\ast,$ 
and $Q_1\sq Q$. 

Assume that $\Ups=X\ti Y$ and $\Ups'=X'\ti Y'.$ Consider the 
non-empty $\is11$ sets $Y_2=\pri P^\ast\sq Y$ and 
$X_2=\ans{x\in X:\sus y\in Y_2\:x\eqf y}.$ It follows from 
Proposition~\ref p that 
$Q_1=\ans{\ang{x,y}\in Q:x\in X_2}\not=\emps,$ hence $Q_1$ is a 
condition in $\dpm$ and $X_1=\pri Q_1$ is a non-empty $\is11$ 
subset of $X_2\sq X$. 

The set $Y_1=\ans{y\in Y_2:\sus x\in X_1\: x\eqf y}$ satisfies 
$X_1\eqf Y_1,$ therefore $\Ups_1=X_1\ti Y_1\in\dpt.$ Furthermore 
$P_1=\ans{\ang{y,x}\in P^\ast:y\in Y_1}\in\dpp$. 

Put $X'_1=\prt P_1\sq X'$ and $Y'_1=\prt Q_1\sq Y'.$ Take notice 
that $Y_1\eqf X'_1$ because any condition in $\dpp$ is a subset 
of $\eqf,$ similarly $X_1\eqf Y'_1,$ and $X_1\eqf Y_1,$ see above. 
It follows that $X'_1\eqf Y'_1,$ hence $\Ups'_1=X'_1\ti Y'_1$ 
is a condition in $\dpt$.

Now $\gp_1=\ang{\Ups_1,P_1,\Ups'_1,Q_1}\in\gP$ as 
required. 

We conclude that $\gP$ forces ``quadruples'' of reals 
$\ang{x,y,x',y'}$ such that the pairs $\ang{x,y}$ and 
$\ang{x',y'}$ are \dd\dpt generic, hence satisfy $x\mek y$ and 
$x'\mek y'$ provided the generic set contains $\Ups_0$ --- by 
the assumption above. Furthermore the pair $\ang{y,x'}$ is 
\dd\dpp generic, 
hence $y\mek x',$ while the pair $\ang{x,y'}$ is \dd\dpm generic, 
hence $x\not\mek y',$ which is a contradiction.
\epf

\parf{The splitting construction}
\label s

Let, in the universe $\kV,$ $\kpa=\gc.$ Let $\kvp$ be a 
\dd\kpa collapse extension of $\kV.$ 

Our aim is to define, in $\kvp,$ a splitting system of sets which 
leads to a function $F$ satisfying \ref{2B} of Theorem~\ref{hb}. 
Let us fix two points before the construction starts. 

{\it First\/}, 
as the forcing notions involved are countable in $\kV,$ there 
exist, in $\kvp,$ enumerations $\ans{D(n):n\in\om},\msur$
$\ans{D_2(n):n\in\om},$ and $\ans{D^2(n):n\in\om},$ of all open 
dense sets in resp.\ $\dP,\msur$ $\dpp,\msur$ $\dpt,$ which (the 
dense sets) belong to $\kV,$ such that $D(n+1)\sq D(n)$ \etc\ for 
each $n$.

{\it Second\/}, we introduce the notion of a crucial pair. 
A pair $\ang{u,v}$ of binary sequences $u,\,v\in 2^n$ is called 
{\it crucial\/} iff $u=1^k\we 0\we w$ and $v=0^k\we 1\we w$ for 
some $k<n$ and $w\in 2^{n-k-1}.$ One easily sees that the graph 
of all crucial pairs in $2^n$ is actually a chain connecting all 
members of $2^n.$ 

We define, in $\kvp,$ a system of sets $X_u\in\dP,$ where 
$u\in 2\lom,$ and sets $P_{uv}\in \dpp,\msur$ $\ang{u,v}$ being a 
crucial pair in some $2^n,$ satisfying the following conditions:

\ben
\itemsep=1mm
\def\theenumi{(\arabic{enumi})}
\def\labelenumi{\theenumi}
\itla1\msur
$X_u\in D(n)$ whenever $u\in 2^n\,;$ $X_{u\we i}\sq X_u$;

\itla2
if $\ang{u,v}$ is a crucial pair in $2^n$ then 
$P_{uv}\in D_2(n)\,;$ $P_{u\we i\,,\,v\we i}\sq P_{uv}$;

\itla3
if $u,\,v\in 2^n$ and $u(n-1)\not=v(n-1)$ then 
$X_u\ti X_v\in\dpt,$ moreover $\in D^2(n)$ and 
$X_u\cap X_v=\emps$;

\itla4
if $\ang{u,v}$ is a crucial pair in $2^n$ then $\pri P_{uv}=X_u$ 
and $\prt P_{uv}=X_v$.
\een

\punk*{\bfit Why this implies the existence of a required 
function ?}


First of all for any $a\in 2^\om$ (in $\kvp$) the sequence 
of sets $X_{a\res n}$ is \dd\dP generic over $\kV$ by \ref1, 
therefore the intersection $\bigcap_{n\in\om}X_{a\res n}$ is 
a singleton. Let $F(a)\in\cN$ be its only element. 

It does not take much effort to prove that $F$ is continuous 
and $1-1$. 

Consider $a,\,b\in 2^\om$ satisfying $a\nEo b.$ Then 
$a(n)\not=b(n)$ for infinitely many $n,$ hence the pair 
$\ang{F(a),F(b)}$ is \dd\dpt generic by \ref3, thus $F(a)$ and 
$F(b)$ are \dd\mek incomparable by Lemma~\ref{pp}.

Consider $a,\,b\in 2^\om$ satisfying $a\meo b.$ We may assume that 
$a$ and $b$ are \dd\meo neighbours, \ie\ $a=1^k\we 0\we c$ while 
$b=0^k\we 1\we c$ for some $k\in\om$ and $c\in 2^\om.$ Then by 
\ref2 the sequence of sets $P_{a\res n\,,\,b\res n},\msur$ 
$n>k,$ is \dd\dpp generic, hence it results in a pair of reals 
satisfying $x\mek y.$ However $x=F(a)$ and $y=F(b)$ by \ref4. 

\punk*{\bfit The construction of a splitting system}

{\it We argue in\/} $\kvp$. 

Suppose that the construction has been completed up to a level 
$n,$ and expand it to the next level. From now on $s,\,t$ will 
denote sequences in $2^n$ while $u,\,v$ will denote sequences 
in $2^{n+1}.$ 

To start with we set $X_{s\we i}=X_s$ for all $s\in 2^n$ and 
$i=0,1,$ and $P_{s\we i\,,\,t\we i}=P_{st}$ whenever $i=0,1$ and 
$\ang{s,t}$ is a crucial pair in $2^n.$ 

For the ``initial'' crucial pair $\ang{1^n\we 0,0^n\we 1}$ at 
this level let 
$P_{1^n\we 0\,,\,0^n\we 1}=X_{1^n\we 0}\ti X_{0^n\we 1}=
X_{1^n}\ti X_{0^n}.$ Then $P_{1^n\we 0\,,\,0^n\we 1}\in\dpt.$
\footnote
{\label{ff}\ It easily follows from \ref2 and \ref4 that 
$X_s\eqf X_t$ for all $s,\,t\in 2^n,$ because $s$ and $t$ are 
connected in $2^n$ by a unique chain of crucial pairs.}  

This ends the definition of ``initial values'' at the 
\dd{n\pone}th level. The plan is to gradually ``shrink'' the 
sets in order to fulfill the requirements.\vom 

{\bfit Step 1\/}. 
We take care of item \ref1. Consider an arbitrary 
$u_0=s_0\we i\in 2^{n+1}.$ As $D(n)$ is dense there is a set  
$X'\in D(n),\msur$ $X'\sq X_{u_0}.$ The intension is to take $X'$ 
as the ``new'' $X_{u_0}.$ But this change has to be expanded 
through the chain of crucial pairs, in order to preserve \ref4. 

Thus put $X'_{u_0}=X'.$ Suppose that $X'_u$ has been defined and 
is included in $X_u,$ the old version, for some $u\in 2^{n+1},$ 
and $\ang{u,v}$ is a crucial pair, $v\in 2^{n+1}$ being not yet 
encountered. Define $P'_{uv}=(X'_u\ti\cN)\cap P_{uv}$ and  
$X'_v=\prt P'_{uv}.$ Clearly \ref4 holds for the new sets 
$X'_u,\msur$ $X'_v,$ and $P'_{uv}$. 

The construction describes how the original change from $X_{u_0}$ 
to $X'_{u_0}$ spreads through the chain of crucial pairs in 
$2^{n+1},$ resulting in a system of new sets, $X'_u$ and 
$P'_{uv},$ which satisfy \ref1 for the particular 
$u_0\in 2^{n+1}.$ 
We iterate this construction consecutively for all 
$u_0\in 2^{n+1},$ getting finally a system of sets satisfying 
\ref1 (fully) (and \ref4), which we shall denote by $X_u$ and 
$P_{uv}$ from now on.\vom

{\bfit Step 2\/}. 
We take care of item \ref3. Fix a pair of $u_0$ and $v_0$ 
in $2^{n+1},$ such that $u_0(n)=0$ and $v_0(n)=1.$ By the density 
of $D^2(n),$ there is a set $X'_{u_0}\ti X'_{v_0}\in D^2(n)$  
included in $X_{u_0}\ti X_{v_0}.$ We may assume that 
$X'_{u_0}\cap X'_{v_0}=\emps.$ 
(Indeed it easily follows from Proposition~\ref p, for $P_-,$ that 
there exist reals $x_0\in X_{u_0}$ and $y_0\in X_{v_0}$ satisfying 
$x_0\eqf y_0$ but $x_0\not=y_0,$ say $x_0(k)=0$ while $y_0(k)=1.$ 
Define 
\dm
X=\ans{x\in X_0:x(k)=0\cj\sus y\in Y_0\:(y(k)=1\cj x\eqf y)},
\dm 
and $Y$ correspondingly; then $X\eqf Y$ and $X\cap Y=\emps$.)

Spread the change from $X_{u_0}$ to $X'_{u_0}$ and from $X_{v_0}$ 
to $X'_{v_0}$ through the chain of crucial pairs in $2^{n+1},$ by 
the method of Step 1, until the wave of spreading from $u_0$ meets 
the wave of spreading from $u_0$ at the ``meeting'' crucial pair 
$\ang{1^n\we 0,0^n\we 1}.$ This 
leads to a system of sets $X'_u$ and $P'_{uv}$ which satisfy 
\ref3 for the particular pair $\ang{u_0,v_0}$ and still 
satisfy \ref4 possibly except for the  ``meeting'' crucial pair 
$\ang{1^n\we 0,0^n\we 1}$ (for which basically 
$P'_{1^n\we 0\,,\,0^n\we 1}$ is not yet defined for this step). 

Note that Step 1 leaves 
$P_{1^n\we 0\,,\,0^n\we 1}$ in the form 
$X_{1^n\we 0}\ti X_{0^n\we 1}$ 
(where $X_{1^n\we 0}$ and $X_{0^n\we 1}$ are the ``versions'' at 
the end of Step 1). We now have new sets, $X'_{1^n\we 0}$ and 
$X'_{0^n\we 1},$ included in resp.\ $X_{1^n\we 0}$ and 
$X_{0^n\we 1}$ and satisfying $X'_{0^n\we 0}\eqf X'_{0^n\we 1}$ 
(because we had $X'_{u_0}\eqf X'_{v_0}$ at the beginning of the 
change.) It remains to define 
$P'_{1^n\we 0\,,\,0^n\we 1}=X'_{1^n\we 0}\ti X'_{0^n\we 1}.$ 
This ends the consideration of the pair 
$\ang{u_0,v_0}$.

Applying this construction consecutively for all pairs of 
$u_0\in P_0$ and $v_0\in P_1$ (including the pair 
$\ang{1^n\we 0,0^n\we 1}$) we finally get a system of sets 
satisfying \ref1, \ref3, and \ref4, which will be 
denoted still by $X_u$ and $P_{uv}$.\vom

{\bfit Step 3\/}.
We finally take care of \ref2. Consider a particular crucial pair 
$\ang{u_0,v_0}$ in $2^{n+1}.$ By the density there 
is a set $P'_{u_0,v_0}\in D_2(n),$ $P'_{u_0,v_0}\sq P_{u_0,v_0}.$ 
(In the case when $\ang{u_0,v_0}$ is 
$\ang{1^n\we 0,0^n\we 1}$ we rather apply Proposition~\ref p to 
obtain $P'_{u_0,v_0}$.)

Define $X'_{u_0}=\pri P'_{u_0,v_0}$ and 
$X'_{v_0}=\prt P'_{u_0,v_0}$ and spread this change through the 
chain of crucial pairs in $2^{n+1}.$ 
(Note that $X'_{u_0}\eqf X'_{v_0}$ as sets in $\dpt$ are 
included in $\eqf.$ This keeps $X'_u\eqf X'_v$ 
for all $u,\,v\in 2^{n+1}$ through the spreading.) 
 
Executing this step for all crucial pairs 
in $2^{n+1}$ we finally end the construction, in $\kvp,$ of a 
system of sets satisfying \ref1 through \ref4.\vtm

\hfill $\squarE$ (Theorem~\ref{hb})


\let\section=\subsection

\end{document}